\documentclass{fpsac}
\usepackage{graphicx}
\usepackage{amssymb}
\usepackage{amsfonts}
\usepackage{epsfig}
\usepackage{amscd}

\newtheorem{example}{Example}[section]

\newtheorem{theorem}[example]{Theorem}

\newtheorem{conjecture}[example]{Conjecture}

\newtheorem{proposition}[example]{Proposition}

\def\MC{M}             % Base standard de CC.
\def\SC{S}             % Algo qui associe une CC a une permutation.

\def\SW{{\mathcal S}}  % Base * des $S$ sur TD.
\def\EW{{\mathcal E}}  % Base * des $E$ sur TD.
\def\MM{{\mathcal M}}  % Base standard de TD.
\def\TD{{\mathfrak{TD}}}% Trigebre dendriforme libre sur un generateur
\def\TC{{\mathfrak{TC}}}% Tricubique libre sur un generateur
\def\TT{{\mathcal T}} % Arbre plan associe a un mot
\def\M{{\bf M}}       % Base standard de WQSym.

\def\pluspetit{\preceq}% Ordre sur le treillis du pseudo-permutoedre.

\def\gaudend{\ll}      % gauche dendriforme
\def\droitdend{\gg}    % droite dendriforme
\def\gautrid{\!\prec\!}   % gauche trigebre dendriforme
\def\miltrid{\circ}       % milieu trigebre dendriforme
\def\droittrid{\!\succ\!} % droite trigebre dendriforme

\def\park{{\bf a}}    % une fonction de parking
\def\Ig{{\bf I}}      % I gras

\def\sep{\,|\,}       % separateurs de l'algebre 3^{n-1} et des parking hypo
\def\limproj{{\rm proj\,lim}}
\def\K{{\mathbb K}}   % corps de base
\def\Sym{{\bf Sym}}   % NCSF

\def\QSym{{\it QSym}}          % QSym
\def\NCSF{{\bf Sym}}           % NCSF
\def\FQSym{{\bf FQSym}}        % permutations
\def\MQSym{{\bf MQSym}}        % matrices
\def\PQSym{{\bf PQSym}}        % Parking 
\def\SQSym{{\bf SQSym}}        % Schroder
\def\SCQSym{{\bf SCQSym}}      % Algebre des 3^n (quotient de SQSym^*)
\def\WQSym{{\bf WQSym}}        % Mots initiaux
\def\NCQSym{{\bf NCQSym}}      % Mots initiaux (pour Bergeron-Zabrocki)
\def\a{{\bf a}}
\def\b{{\bf b}}

\def\pack{{\rm pack}}          % packed word
\def\mup{{\rm mup}}            % maximally unpacked word
\def\last{{\rm last}}   % last letter
\def\ssh{\Cup}          % shuffle shifte
\def\sconc{\bullet}     % concatenation shiftee
\def\Std{{\rm Std}}     % standardisation
\def\Park{{\rm Park}}   % parkisation
\def\convW{{*_W}}       % une loi de convolution
\def\<{\langle}
\def\>{\rangle}
\def\park{{\bf a}} % les parking 
\def\F{{\bf F}}         % F de FQSym, PQSym, ...
\def\G{{\bf G}}         % G de FQSym^*
\def\SG{{\mathfrak S}}  % groupe symetrique

\def\dim{{\rm dim}}

\def\PW{{\rm PW}}   % packed words = mots tasses

\def\shuff#1#2{\mathbin{
      \hbox{\vbox{
        \hbox{\vrule
              \hskip#2
              \vrule height#1 width 0pt
               }%
        \hrule}%
             \vbox{
        \hbox{\vrule
              \hskip#2
              \vrule height#1 width 0pt
               \vrule }%
        \hrule}%
}}}
\def\shuffl{{\mathchoice{\shuff{7pt}{3.5pt}}%
                        {\shuff{6pt}{3pt}}%
                        {\shuff{4pt}{2pt}}%
                        {\shuff{3pt}{1.5pt}}}}%
\def\shuffle{\, \shuffl \,}

\title{Polynomial realizations of some trialgebras}
\author%[J.-C.~Novelli and J.-Y.~Thibon]
{Jean-Christophe Novelli and Jean-Yves Thibon}

\address[] {Institut Gaspard Monge, Universit\'e de Marne-la-Vall\'ee \\
5 Boulevard Descartes \\Champs-sur-Marne \\77454 Marne-la-Vall\'ee cedex 2 \\
FRANCE}
\email[Jean-Christophe Novelli]{novelli@univ-mlv.fr}
\email[Jean-Yves Thibon]{jyt@univ-mlv.fr} 
\date{}



\begin{document}

\begin{abstract}
We realize several combinatorial Hopf algebras based on set compositions, plane
trees and segmented compositions in terms of noncommutative polynomials in
infinitely many variables. For each of them, we describe a trialgebra
structure, an internal product, and several bases.

\begin{resume}
Nous r\'ealisons plusieurs alg\`ebres de Hopf combinatoires dont les bases
sont index\'ees par les partitions d'ensembles ordonn\'ees, les arbres plans
et les compositions segment\'ees en termes de polyn\^omes non-commutatifs en
une infinit\'e de variables. Pour chacune d'elles, nous d\'ecrivons sa
structure de trig\`ebre, un produit int\'erieur et plusieurs bases.
\end{resume}
\end{abstract}

\maketitle

%%%%%%%%%%%%%%%%%%%%%%%%%%%%%%%%%%%%%%%%%%%%%%%%%%%%%%%%%%%%%%%%%%%%%%%%%%%%%%%
%%%%%%%%%%%%%%%%%%%%%%%%%%%%%%%%%%%%%%%%%%%%%%%%%%%%%%%%%%%%%%%%%%%%%%%%%%%%%%%
%%%%%%%%%%%%%%%%%%%%%%%%%%%%%%%%%%%%%%%%%%%%%%%%%%%%%%%%%%%%%%%%%%%%%%%%%%%%%%%
\section{Introduction}

The aim of this note is to construct and analyze several combinatorial
Hopf algebras arising in the theory of operads from the point of view of the
theory of noncommutative symmetric functions.
Our starting point will be the algebra of noncommutative polynomial invariants
$$\WQSym(A)=\K\<A\>^{\SG(A)_{QS}}$$ of Hivert's quasi-symmetrizing action
\cite{Hiv}.
It is known that, when the alphabet $A$ is infinite, $\WQSym(A)$ acquires the
structure of a graded Hopf algebra whose bases are parametrized by ordered set
partitions (also called set compositions) \cite{Hiv,NT4,BZ}. Set compositions
are in one-to-one correspondence with faces of permutohedra, and actually,
$\WQSym$ turns out to be isomorphic to one of the Hopf algebras introduced by
Chapoton in \cite{Chap}.  From this algebra, Chapoton obtained graded Hopf
algebras based on the faces of the associahedra (corresponding to plane trees
counted by the little Schr\"oder numbers) and on faces of the hypercubes
(counted by powers of 3).
Since then, Loday and Ronco have introduced the operads of dendriform
trialgebras and of tricubical algebras \cite{LRtri}, in which the free
algebras on one generator are respectively based on faces of associahedras and
hypercubes, and are isomorphic (as Hopf algebras) to the corresponding
algebras of Chapoton.  More recently, we have introduced a Hopf algebra
$\PQSym$, based on parking functions \cite{NT1,NT2,NT3}, and derived from it a
series of Hopf subalgebras or quotients, some of which being isomorphic to the
above mentioned ones as associative algebras, but not as Hopf algebras.

In the following, we will show that applying the same techniques,
starting from $\WQSym$ instead of $\PQSym$, allows one to recover all of these
algebras, together  with their original Hopf structure, in a very natural way.
This provides in particular for each of them an explicit realization in terms
of noncommutative polynomials.  The Hopf structures can be analyzed very
efficiently by means of Foissy's theory of bidendriform bialgebras~\cite{Foi}.
A natural embedding of $\WQSym$ in $\PQSym^*$ implies that $\WQSym$ is
bidendriform, hence, free and self-dual. These properties are inherited by
$\TD$, the free dendriform trialgebra on one generator, and some of them by
$\TC$, the free cubical trialgebra on one generator.
A lattice structure on the set of faces of the permutohedron (introduced in
\cite{KLNPS} under the name  ``pseudo-permutohedron'' and rediscovered
in~\cite{PR}) leads to the construction of various bases of these algebras.
Finally, the natural identification of the homogeneous components of the dual
$\WQSym^*_n$ (endowed with the internal product induced by $\PQSym$) with the
Solomon-Tits algebras (that is, the face algebras of the braid arrangements of
hyperplanes) implies that all three algebras admit an internal product.

{\footnotesize
\bigskip {\it Notations -- } We  assume that the reader is familiar with
the standard notations of the theory of noncommutative symmetric functions
\cite{NCSF1,NCSF6} and with the Hopf algebra of parking functions
\cite{NT1,NT2,NT3}.
%%%
We shall need an infinite totally ordered alphabet
$A=\{a_1<a_2<\cdots<a_n<\cdots\}$, generally assumed to be the set of
positive integers.
We denote by $\K$ a field of characteristic $0$, and by $\K\<A\>$ the free
associative algebra over $A$ when $A$ is finite, and the projective limit
$\limproj_B \K\<B\>$, where $B$ runs over finite subsets of $A$, when $A$ is
infinite.
The \emph{evaluation} of a word $w$ is the sequence whose $i$-th term is the
number of times the letter $a_i$ occurs in $w$.
The \emph{standardized word} $\Std(w)$ of a word $w\in A^*$ is the permutation
obtained by iteratively scanning $w$ from left to right, and labelling
$1,2,\ldots$ the occurrences of its smallest letter, then numbering the
occurrences of the next one, and so on.
% Alternatively, $\sigma=\Std(w)^{-1}$
%can be characterized as the unique permutation of minimal length such that
%$w\sigma$ is a nondecreasing word.
For example, $\Std(bbacab)=341624$.
For a word $w$ on the alphabet $\{1,2,\ldots\}$, we denote by $w[k]$ the word
obtained by replacing each letter $i$ by the integer $i+k$.
If $u$ and $v$ are two words, with $u$ of length $k$, one defines
the \emph{shifted concatenation}
$u\sconc v = u\cdot (v[k])$
and the \emph{shifted shuffle}
$ u\ssh v= u\shuffle (v[k])$,
where $\shuffle$ is the usual shuffle product.
}

%%%%%%%%%%%%%%%%%%%%%%%%%%%%%%%%%%%%%%%%%%%%%%%%%%%%%%%%%%%%%%%%%%%%%%%%%%%%%%%
%%%%%%%%%%%%%%%%%%%%%%%%%%%%%%%%%%%%%%%%%%%%%%%%%%%%%%%%%%%%%%%%%%%%%%%%%%%%%%%
%%%%%%%%%%%%%%%%%%%%%%%%%%%%%%%%%%%%%%%%%%%%%%%%%%%%%%%%%%%%%%%%%%%%%%%%%%%%%%%
\section{The Hopf algebra $\WQSym$}

%%%%%%%%%%%%%%%%%%%%%%%%%%%%%%%%%%%%%%%%%%%%%%%%%%%%%%%%%%%%%%%%%%%%%%%%%%%%%%%
\subsection{Noncommutative quasi-symmetric invariants}

The \emph{packed word} $u=\pack(w)$ associated with a word $w\in A^*$ is
obtained by the following process. If $b_1<b_2<\ldots <b_r$ are the letters
occuring in $w$, $u$ is the image of $w$ by the homomorphism
$b_i\mapsto a_i$.
A word $u$ is said to be \emph{packed} if $\pack(u)=u$. We denote by $\PW$ the
set of packed words.
With such a word, we associate the polynomial
\begin{equation}
\M_u :=\sum_{\pack(w)=u}w\,.
\end{equation}
For example, restricting $A$ to the first five integers,
\begin{equation}
%\begin{split}
\M_{13132}= 13132 + 14142 + 14143 + 24243 
+ 15152 + 15153 + 25253 + 15154 + 25254 + 35354.
%\end{split}
\end{equation}
Under the abelianization
$\chi:\ \K\langle A\rangle\rightarrow\K[X]$, the $\M_u$ are mapped to the
monomial quasi-symmetric functions $M_I$ 
($I=(|u|_a)_{a\in A}$ being the evaluation vector of $u$).

These polynomials span a subalgebra 
of $\K\langle A\rangle$, called $\WQSym$ for Word
Quasi-Symmetric functions~\cite{Hiv} (and called $\NCQSym$ in~\cite{BZ}),
consisting in the invariants of the noncommutative
version of Hivert's quasi-symmetrizing action \cite{Hiv-adv}, which is defined
by
$\sigma\cdot w = w'$ where $w'$ is such that $\Std(w')=\Std(w)$ and
$\chi(w')=\sigma\cdot\chi(w)$.
Hence, two words are in the same
$\SG(A)$-orbit iff they have the same packed word.

$\WQSym$ can be embedded in $\MQSym$~\cite{Hiv,NCSF6}, by
$\M_u\mapsto {\bf MS}_M$,
where $M$ is the packed $(0,1)$-matrix whose $j$th column contains exactly
one 1 at row $i$ whenever the $j$th letter of $u$ is $a_i$. Since the duality
in $\MQSym$ consists in tranposing the matrices, one can also embed $\WQSym^*$
in $\MQSym$.
The multiplication formula for the basis $\M_u$ follows from
that of ${\bf MS}_M$ in $\MQSym$: 
\begin{proposition}
The product on $\WQSym$ is given by
\begin{equation} 
\label{prodG-wq}
\M_{u'} \M_{u''} = \sum_{u \in u'\convW u''} \M_u\,,
\end{equation}
where the \emph{convolution} $u'\convW u''$ of two packed words
is defined as
\begin{equation} 
u'\convW u'' = \sum_{v,w ;
u=v\cdot w\,\in\,\PW, \pack(v)=u', \pack(w)=u''} u\,.
\end{equation}
\end{proposition}

For example,
\begin{equation}
\M_{11} \M_{21} =
\M_{1121} + \M_{1132} + \M_{2221} + \M_{2231} + \M_{3321}.
\end{equation}
%\begin{equation}
%\begin{split}
%\M_{21} \M_{121} =& \ \ \
%\M_{12121} + \M_{12131} + \M_{12232} + \M_{12343} + \M_{13121} +
%\M_{13232} + \M_{13242}\\
%&+ \M_{14232} + \M_{23121} + \M_{23131} + \M_{23141} + \M_{24131} +
%\M_{34121}.
%\end{split}
%\end{equation}
%
Similarly, the embedding in $\MQSym$ implies immediately
 that $\WQSym$ is a Hopf subalgebra of
$\MQSym$. However, the coproduct can also be defined directly by the usual
trick of noncommutative symmetric functions, considering the
alphabet $A$ as an ordered sum of two mutually commuting alphabets 
$A'\hat+A''$. First, by direct inspection, one finds that 
\begin{equation}
\M_u(A'\hat+A'') = \sum_{0\leq k\leq \max(u)}
\M_{(u|_{[1,k]})}(A') \M_{\pack(u|_{[k+1,\max(u)})}(A''),
\end{equation}
where $u|_{B}$ denote the subword obtained by restricting $u$ to the
subset $B$ of the alphabet,
and now, the coproduct $\Delta$ defined by
\begin{equation}
\Delta \M_u(A)
=\sum_{0\leq k\leq \max(u)}
\M_{(u|_{[1,k]})}\otimes \M_{\pack(u|_{[k+1,\max(u)})},
\end{equation}
is then clearly a morphism for the concatenation product, hence defines a bialgebra
structure.

Given two packed words $u$ and $v$, define the \emph{packed shifted shuffle}
$u\ssh_W\, v$ as the shuffle product of $u$ and $v[\max(u)]$.
One then easily sees that
\begin{equation}
\Delta \M_w(A)= \sum_{u,v ; w\in u\ssh_W\,v} \M_u\otimes \M_v.
\end{equation}
For example,
%\begin{equation}
%\Delta \M_{121} = 1\otimes\M_{121} + \M_{11}\otimes\M_{1} + \M_{121}\otimes1.
%\end{equation}
%
\begin{equation}
\label{D-32121}
\Delta \M_{32121} = 1\otimes\M_{32121} + \M_{11}\otimes\M_{211} +
\M_{2121}\otimes\M_{1} + \M_{32121}\otimes1.
\end{equation}
%
%\begin{equation}
%\begin{split}
%\Delta \M_{1432152} =&\ \ 
%1\otimes \M_{1432152} + \M_{11}\otimes\M_{32141} + \M_{1212}\otimes\M_{213}\\
%&+ \M_{13122}\otimes\M_{12} + \M_{143212}\otimes\M_{1} + \M_{1432152}\otimes1.
%\end{split}
%\end{equation}

Packed words can be naturally identified with \emph{ordered set partitions},
the letter $a_i$ at the $j$th position meaning that $j$ belongs to block $i$.
For example,
\begin{equation}
%\label{init-setpart}
u=313144132 \ \leftrightarrow\ \Pi=(\{2,4,7\},\{9\},\{1,3,8\},\{5,6\})\,.
\end{equation}
To improve the readability of the formulas, we write instead of $\Pi$ a
\emph{segmented permutation}, that is, the permutation obtained by reading the
blocks of $\Pi$ in increasing order and inserting bars $|$ between blocks.

For example,
\begin{equation}
\label{init-setpart}
\Pi=(\{2,4,7\},\{9\},\{1,3,8\},\{5,6\})\,\
\leftrightarrow 247|9|138|56.
\end{equation}
On this representation, the coproduct amounts to deconcatenate the blocks, and
then standardize the factors.
For example, in terms of segmented permutations, Equation~(\ref{D-32121})
reads
\begin{equation}
%\begin{split}
\Delta \M_{35|24|1} = %& \ \ \ \
1\otimes\M_{35|24|1}
+ \M_{12}\otimes\M_{23|1}%\\
+\M_{24|13}\otimes\M_{1}
+ \M_{35|24|1}\otimes1.
%\end{split}
\end{equation}

The dimensions of the homogeneous components of $\WQSym$ are the ordered
Bell numbers $1$, $1$, $3$, $13$, $75$, $541$, $\ldots$
(sequence A000670,~\cite{Slo})
%\begin{equation}
%1 + t + 3\,t^2 + 13\,t^3 + 75\,t^4 + 541\,t^5 + 4\,683\,t^6 + 47\,293\,t^7
%+ 545\,835\,t^8 + 7\,087\,261\,t^9 + \ldots
%\end{equation}
%
so that
\begin{equation}
\dim \WQSym_n = \sum_{k=1}^nS(n,k)k! = A_n(2)\,,
\end{equation}
where $A_n(q)$ are the Eulerian polynomials.

%%%%%%%%%%%%%%%%%%%%%%%%%%%%%%%%%%%%%%%%%%%%%%%%%%%%%%%%%%%%%%%%%%%%%%%%%%%%%%%
\subsection{The trialgebra structure of $\WQSym$} \label{trigd1}

A \emph{dendriform trialgebra} \cite{LRtri} is an associative algebra whose multiplication
$\odot$ splits into three pieces
\begin{equation}
x\odot y = x\gautrid y + x\miltrid y + x\droittrid y\,,
\end{equation}
where $\miltrid$ is associative, and
%\begin{eqnarray}
%(x\gautrid y)\gautrid z = x\gautrid (y\odot z)\,,\\
%(x\droittrid y)\gautrid z = x\droittrid (y\gautrid z)\,,\\
%(x\odot y)\droittrid z = x\droittrid (y\droittrid z)\,,\\
%(x\droittrid y)\miltrid z = x\droittrid (y\miltrid z)\,,\\
%(x\gautrid y)\miltrid z = x\miltrid (y\droittrid z)\,,\\
%(x\miltrid y)\gautrid z = x\miltrid (y\gautrid z)\,.
%\end{eqnarray}

\begin{equation}
(x\gautrid y)\gautrid z = x\gautrid (y\odot z)\,,\ \
(x\droittrid y)\gautrid z = x\droittrid (y\gautrid z)\,,\ \
(x\odot y)\droittrid z = x\droittrid (y\droittrid z)\,,\ \
\end{equation}
\begin{equation}
(x\droittrid y)\miltrid z = x\droittrid (y\miltrid z)\,,\ \
(x\gautrid y)\miltrid z = x\miltrid (y\droittrid z)\,,\ \
(x\miltrid y)\gautrid z = x\miltrid (y\gautrid z)\,.
\end{equation}

%Let $A=\{a_1<a_2<\cdots <a_n <\cdots \,\}$ be an infinite linearly ordered
%alphabet. Here $\K\langle A\rangle$ is understood as the projective
%limit of the $\K\langle A_n\rangle$ where $A_n$ is the interval $[a_1,a_n]$ of
%$A$.
%We denote by $\max(w)$ the greatest letter occuring in the word $w\in A^*$.

It has been shown in~\cite{NT3} that the augmentation ideal
$\K\langle A_n\rangle^+$ has a natural structure of dendriform trialgebra: 
%\medskip
%\begin{definition}[\cite{NT3}]
%\label{def-trigdend}
for two non empty words $u,v\in A^*$, we set
\begin{eqnarray}
u\gautrid v=\begin{cases} uv &
\text{if $\max(u)>\max(v)$}\cr 0 &\mbox{otherwise,} \end{cases}\\
u\miltrid v=\begin{cases} uv &
\text{if $\max(u)=\max(v)$}\cr 0 &\mbox{otherwise,} \end{cases}\\
u\droittrid v=\begin{cases} uv &
\text{if $\max(u)<\max(v)$}\cr 0 &\mbox{otherwise.} \end{cases}
\end{eqnarray}
%\end{definition}

%\begin{lemma}[\cite{NT3}]
%\label{lem-trigdend}
%The three operations $\gautrid\,, \miltrid,\droittrid\,$, endow the
%augmentation ideal $\K\<A\>^+$ with the structure of a dendriform trialgebra.
%\end{lemma}
%%%
\begin{theorem}
\label{wqs-trid}
$\WQSym^+$ is a sub-dendriform trialgebra of $\K\<A\>^+$,
the partial products being given by
\begin{equation}
\M_{w'} \gautrid \M_{w''} =
\sum_{w=u.v\in w'\convW w'', |u|=|w'| ; \max(v)<\max(u)}
\M_w,
\end{equation}
\begin{equation}
\M_{w'} \miltrid \M_{w''} =
\sum_{w=u.v\in w'\convW w'', |u|=|w'| ; \max(v)=\max(u)}
\M_w,
\end{equation}
\begin{equation}
\M_{w'} \droittrid \M_{w''} =
\sum_{w=u.v\in w'\convW w'', |u|=|w'| ; \max(v)>\max(u)}
\M_w,
\end{equation}
\end{theorem}
It is known~\cite{LRtri} that the free dendriform trialgebra on one
generator, denoted here by
$\TD$, is a free associative algebra with Hilbert series
\begin{equation}
\label{sg-dendt}
\sum_{n\geq0} s_n t^n = \frac{1+t-\sqrt{1-6t+t^2}}{4t}
= 1 + t + 3t^2 + 11t^3 + 45t^4 + 197t^5 + \cdots,
\end{equation}
the generating function of the \emph{super-Catalan}, or
\emph{little Schr\"oder} numbers, counting \emph{plane trees}.
The previous considerations allow us to give a simple polynomial realization
of $\TD$. Consider the polynomial
\begin{equation}
\M_1=\sum_{i\ge 1}a_i \ \in \WQSym \,,
\end{equation}
%(the sum of all letters). We can then state:
%
\begin{theorem}[\cite{NT3}]
The sub-trialgebra $\TD$ of $\WQSym^+$ generated by $\M_1$ is free as a
dendriform trialgebra.
\end{theorem}
%
%We will see in Section~\ref{sqsym}  other structures 
%on $\TD$ implied by this realization.
%
Based on numerical evidence, we conjecture the following result:

\begin{conjecture}
$\WQSym$ is a free dendriform trialgebra.
\end{conjecture}
%
%Applying the same trick as in~\cite{Foi} for computing the generating
%series of the totally primitive elements, one gets the generating series of
The number $g'_n$ of generators in degree $n$ of $\WQSym$ as a free
dendriform trialgebra would then be
\begin{equation}
%\begin{split}
\sum_{n\geq0} g'_nt^n = \frac{OB(t)-1}{2OB(t)^2 - OB(t)}% \\
= t + 2\,t^3 + 18\, t^4 + 170\,t^5 + 1\,794\,t^6 + 21\,082\,t^7 +O(t^8).%\\
%&+ 273\,714\,t^8 + 3\,897\,994\,t^9 + O(t^{10}).
%\end{split}
\end{equation}
where $OB(t)$ is the generating series of the ordered Bell numbers.

%%%%%%%%%%%%%%%%%%%%%%%%%%%%%%%%%%%%%%%%%%%%%%%%%%%%%%%%%%%%%%%%%%%%%%%%%%%%%%%
\subsection{Bidendriform structure of $\WQSym$}

A \emph{dendriform dialgebra}, as defined by Loday~\cite{Lod-dend}, is an
associative algebra $D$ whose multiplication $\odot$ splits into two binary
operations
\begin{equation}
x \odot y = x\gaudend y + x\droitdend y\,,
\end{equation}
called left and right, satisfying the  following three compatibility relations
for all $a$, $b$, and $c$ different from $1$ in $D$:
\begin{equation}
\label{dend1}
(a\gaudend b)\gaudend c = a\gaudend (b\odot c),
\quad
%\end{equation}
%\begin{equation}
%\label{dend2}
(a\droitdend b)\gaudend  c = a\droitdend  (b\gaudend c),
%\end{equation}
%\begin{equation}
%\label{dend3}
\quad
(a\odot b)\droitdend  c = a\droitdend  (b\droitdend c).
\end{equation}

%%
%A \emph{bidendriform bialgebra} \cite{Foi} is a  dendriform dialgebra equipped
%with a coproduct which splits into two parts, satifying the
%\emph{codendriform relations}, obtained by dualizing the dendriform
%relations, and certain compatibility properties with the two half-products.
%
A \emph{codendriform coalgebra} is a coalgebra $C$ whose coproduct $\Delta$
splits as $\Delta(c)=\overline\Delta(c)+c\otimes1+1\otimes c$ and
$\overline\Delta=\Delta_\gaudend  + \Delta_\droitdend $, such that, for all
$c$ in $C$:
\begin{equation}
\label{codend1}
(\Delta_\gaudend \otimes Id) \circ \Delta_\gaudend (c) =
(Id\otimes\overline\Delta)\circ\Delta_\gaudend (c),
\end{equation}
\begin{equation}
\label{codend2}
(\Delta_\droitdend \otimes Id) \circ \Delta_\gaudend (c) =
(Id\otimes\Delta_\gaudend )\circ\Delta_\droitdend (c),
\end{equation}
\begin{equation}
\label{codend3}
(\overline\Delta\otimes Id)\circ\Delta_\droitdend (c) =
(Id\otimes\Delta_\droitdend ) \circ \Delta_\droitdend (c).
\end{equation}

The Loday-Ronco algebra of planar binary trees introduced in~\cite{LR1} arises
as the free dendriform dialgebra on one generator. This is moreover a Hopf
algebra, which turns out to be self-dual, so that it is also codendriform.
There is some compatibility between the dendriform and the codendriform
structures, leading to what has been called by Foissy~\cite{Foi} a
\emph{bidendriform bialgebra}, defined as a bialgebra which is
both a dendriform dialgebra and a codendriform
coalgebra, satisfying the following four compatibility relations
\begin{equation}
\label{bidend1}
\Delta_\droitdend  (a\droitdend b) =
      a'b'_\droitdend \!\otimes\!  a''\!\!\droitdend \!b''_\droitdend \,
\,+\, a'\!\otimes\! a''\!\!\droitdend \!b \,
\,+\, b'_\droitdend \!\otimes\! a\!\droitdend \!b''_\droitdend  \,
\,+\, ab'_\droitdend \!\otimes\! b''_\droitdend  \,
\,+\, a\!\otimes\! b\,,
\end{equation}
\begin{equation}
\label{bidend2}
\Delta_\droitdend  (a\gaudend b) =
      a'b'_\droitdend \!\otimes\! a''\!\gaudend \!b''_\droitdend 
\,+\, a'\!\otimes\! a''\!\gaudend \!b 
\,+\, b'_\droitdend \!\otimes\! a\!\gaudend \!b''_\droitdend \,,
\end{equation}
\begin{equation}
\label{bidend3}
\Delta_\gaudend  (a\droitdend b) =
      a'b'_\gaudend \!\otimes a''\!\droitdend \!b''_\gaudend
\,+\, ab'_\gaudend \otimes b''_\gaudend  
\,+\, b'_\gaudend \otimes a\droitdend b''_\gaudend \,,
\end{equation}
\begin{equation}
\label{bidend4}
\Delta_\gaudend  (a\gaudend b) =  a'b'_\gaudend \!\otimes\! a''\!\gaudend
\!b''_\gaudend 
\,+\, a'b\!\otimes\! a'' 
\,+\, b'_\gaudend \!\otimes\! a\!\gaudend \!b''_\gaudend  \,+\, b\!\otimes\! a\,,
\end{equation}
where the pairs $(x',x'')$ (resp. $(x'_\gaudend ,x''_\gaudend )$ and
$(x'_\droitdend,x''_\droitdend)$) correspond to all possible elements
occuring in $\overline\Delta x$ (resp. $\Delta_\gaudend x$ and
$\Delta_\droitdend  x$), summation signs being understood (Sweedler's
notation).

Foissy has shown~\cite{Foi} that a connected bidendriform bialgebra
${\mathcal B}$ is always free as an associative algebra and self-dual as a
Hopf algebra. Moreover, its primitive Lie algebra is free, and as a dendriform
dialgebra, ${\mathcal B}$ is also free over the space of totally primitive
elements (those annihilated by $\Delta_\gaudend$ and $\Delta_\droitdend$).
It is also proved in \cite{Foi} that $\FQSym$ is bidendriform, so that it
satisfies all these properties. In \cite{NT3}, we have proved that
$\PQSym$, the Hopf algebra of parking functions, as also bidendriform.

The realization of $\PQSym^*$ given in \cite{NT2,NT3}  implies that
\begin{equation}
\M_u=\sum_{\pack(\a)=u}\G_\a\,.
\end{equation}
Hence, $\WQSym$ is a subalgebra of $\PQSym^*$. Since in both cases the
coproduct correponds to $A\rightarrow A' \hat + A''$, it is actually a Hopf
subalgebra. It also stable by the
tridendriform operations, and by the codendriform half-coproducts. Hence,

\begin{theorem}
\label{wqs-bid}
$\WQSym$ is a  sub-bidendriform
bialgebra of $\PQSym^*$. More precisely, the product rules are
\begin{equation}
\M_{w'} \gaudend \M_{w''} =
\sum_{w=u.v\in w'\convW w'', |u|=|w'| ; \max(v)<\max(u)}
\M_w,
\end{equation}
\begin{equation}
\M_{w'} \droitdend \M_{w''} =
\sum_{w=u.v\in w'\convW w'', |u|=|w'| ; \max(v)\geq\max(u)}
\M_w,
\end{equation}
\begin{equation}
\Delta_\gaudend \M_w = \sum_{w\in u\ssh_W\, v ; \last(w)\leq |u|}
\M_u\otimes \M_v,
\end{equation}
\begin{equation}
\Delta_\droitdend \M_w = \sum_{w\in u\ssh_W\, v ; \last(w)>|u|}
\M_u\otimes \M_v.
\end{equation}
where $|u|\geq1$ and $|v|\geq1$, and $\last(w)$ means the
last letter of $w$.
As a consequence, $\WQSym$ is free, cofree, self-dual, and its primitive Lie
algebra is free.
%\qed
\end{theorem}

%%%%%%%%%%%%%%%%%%%%%%%%%%%%%%%%%%%%%%%%%%%%%%%%%%%%%%%%%%%%%%%%%%%%%%%%%%%%%%%
\subsection{Duality: embedding $\WQSym^*$ into $\PQSym$}

Recall  from \cite{NT1} that $\PQSym$ is the algebra
with basis $(\F_\park)$,
the product being given by the shifted shuffle of  parking functions,
and that $(\G_\park)$ is the dual basis in $\PQSym^*$.
 
For a packed word $u$ over the integers, let us define its
\emph{maximal unpacking} $\mup(u)$ as the greatest parking function $\b$
for the lexicographic order such that $\pack(\b)=u$.
For example, $\mup(321412451)=641714791$.

Since the basis $(\M_u)$ of $\WQSym$ can be expressed as the sum of $\G_\park$
with a given packed word, the dual basis of $(\M_u)$ in $\WQSym^*$ can be
identified with equivalence classes of  $(\F_\park)$ under the relation
$\F_\park=\F_{\park'}$ iff $\pack(\park)=\pack(\park')$.
Since the shifted shuffle of two maximally unpacked parking functions 
contains only  maximally unpacked parking functions, the dual algebra
$\WQSym^*$ is  in fact a subalgebra of $\PQSym$. Finally, since, if $\park$ is
maximally unpacked then only maximally unpacked parking functions appear in
the coproduct $\Delta\F_\park$, one has

\begin{theorem}
$\WQSym^*$ is a Hopf subalgebra of $\PQSym$. Its  basis element $\M_u^*$ 
can be identified with $\F_\b$ where $\b=\mup(u)$.
\end{theorem}

So we have
\begin{equation}
\label{prodF-WQ}
\F_{\b'}\F_{\b'}:=\sum_{\b\in\b'\ssh\b''}\F_\b\,, \qquad
%\end{equation}
%\begin{equation}
\Delta \F_{\b} =
  \sum_{u\cdot v=\b} \F_{\Park(u)} \otimes \F_{\Park(v)}\,,
\end{equation}
where $\Park$ is the parkization algorithm defined in~\cite{NT3}.
For example,
\begin{equation}
\F_{113} \F_{11} = \F_{11344} + \F_{11434} + \F_{11443} + \F_{14134} 
+ \F_{14143} + \F_{14413} + \F_{41134} + \F_{41143} + \F_{41413} + \F_{44113}.
\end{equation}
\begin{equation}
\Delta\F_{531613} = 1\otimes\F_{531613} + \F_{1}\otimes\F_{31513} +
\F_{21}\otimes\F_{1413} + \F_{321}\otimes\F_{312} + \F_{3214}\otimes\F_{12} +
\F_{43151}\otimes\F_{1} \F_{531613}\otimes1.
\end{equation}

%%%%%%%%%%%%%%%%%%%%%%%%%%%%%%%%%%%%%%%%%%%%%%%%%%%%%%%%%%%%%%%%%%%%%%%%%%%%%%%
\subsection{The Solomon-Tits algebra}

The above realization of $\WQSym^*$ in $\PQSym$ is stable under the internal
product of $\PQSym$ defined in \cite{NT2}.
Indeed, by definition of the internal product, if $\b'$ and $\b''$ are
maximally unpacked, and $\F_\b=\F_{\b'}*\F_{\b''}$, then $\b$ is also
maximally unpacked.

Moreover, if one writes $\b'=\{s'_1,\ldots,s'_k\}$ and $\b''=\{s''_1,\ldots,s''_l\}$ as
ordered set partitions, then the parkized word $\b=\Park(\b',\b'')$ 
corresponds to the ordered set
partition obtained from
\begin{equation}
\{ s'_1\cap s''_1, s'_1\cap s''_2, \ldots, s'_1\cap s''_l,
   s'_2\cap s''_1, \ldots, s'_k\cap s''_l
   \}.
\end{equation}
This formula was rediscovered in~\cite{BZ} and Bergeron and Zabrocki
recognized the Solomon-Tits algebra, in the version given by
Bidigare~\cite{Bid}, in terms of the face semigroup of the braid arrangement
of hyperplanes. So,

\begin{theorem}
$(\WQSym^*,*)$ is isomorphic to the Solomon-Tits algebra.
\end{theorem}
In particular, the product of the Solomon-Tits algebra is dual to the
coproduct $\delta \G(A)=\G(A'A'')$.

%%%%%%%%%%%%%%%%%%%%%%%%%%%%%%%%%%%%%%%%%%%%%%%%%%%%%%%%%%%%%%%%%%%%%%%%%%%%%%%
\subsection{The pseudo-permutohedron}

We shall now make use of the lattice of pseudo-permutations, a combinatorial
structure defined in~\cite{KLNPS} and rediscovered in~\cite{PR}. 
\emph{Pseudo-permutations} are nothing but  ordered set partitions. However,
regarding them as generalized permutations helps uncovering their lattice
structure. Indeed, let us say that if $i$ is in a block strictly to the right
of $j$ with $i<j$ then we have a full inversion $(i,j)$, and
that if $i$ is in the same block as $j$, then we have a \emph{half} inversion
$\frac12 (i,j)$. The total number of inversions is the sum of these numbers.
For example, the table of inversions of $45|13|267|8$ is 
\begin{equation}
\left\{\frac{1}{2}(1,3),\ (1,4),\ (1,5),\ (2,3),\ (2,4),\ (2,5),\
\frac12(2,6),\ \frac12(2,7),\
   (3,4),\ (3,5),\ \frac{1}{2}(4,5),\ \frac12(6,7)\right\},
\end{equation}
and it has $9.5$ inversions.

One can now define a partial order $\pluspetit$ on pseudo-permutations
by setting $p_1\pluspetit p_2$ if the value of the
inversion $(i,j)$ in the table of inversions of $p_1$ is smaller than or equal
to its value in the table of inversions of $p_2$, for all $(i,j)$. This
partial order is a lattice~\cite{KLNPS}.
%To make it compatible with our algebra $\WQSym$, we will consider a slight
%variation of this lattice, that amounts to apply a simple involution to its
%vertices, thus not changing the fact it is a lattice.
%
In terms of packed words, the covering relation 
reads as follows. The successors of a packed word $u$ are the packed
words $v$ such that
\begin{itemize}
\item if all the $i-1$ are to the left of all the $i$ in $u$ then $u$ has as
successor the element where all letters $j$ greater than or equal to $i$ are
replaced by $j-1$.
\item if there are $k$ letters $i$ in $u$, then one can choose an integer $j$
in the interval $[1,k-1]$ and change the $j$ righmost letters $i$  into $i+1$
and the letters $l$ greater than $i$ into $l+1$.
\end{itemize}

For example, $w=44253313$ has five successors,
\begin{equation}
33242212,\ 44243313,\ 55264313,\ 55264413,\ 54263313.
\end{equation}
%The lattice for $n=3$ appears on Figure~\ref{pseudo3}.

\begin{figure}[ht]
\resizebox{3cm}{4cm}{\includegraphics{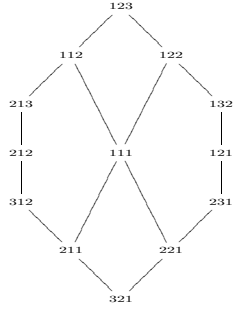}}
%\centerline{\epsfxsize=6cm \epsfbox{pseudo3.eps}}
\caption{\label{pseudo3}The pseudo-permutohedron of degree $3$.}
\end{figure}
%
%The first relation between this lattice and $\WQSym$ is
%
\begin{theorem}[\cite{PR}]
Let $u$ and $v$ be two packed words.
Then $\M_u\M_v$ is an interval of the pseudo-permutohedron lattice.
The minimum of the interval is given by $u\cdot  v[\max(u)]$ and its maximum by
$u[\max(v)]\cdot v$.
\end{theorem}

For example, 
\begin{equation}
\M_{13214} \M_{212} = \sum_{u\in [ 13214656, 35436212 ]}\M_u.
\end{equation}

%%%%%%%%%%%%%%%%%%%%%%%%%%%%%%%%%%%%%%%%%%%%%%%%%%%%%%%%%%%%%%%%%%%%%%%%%%%%%%%
\subsection{Other bases of $\WQSym$ and $\WQSym^*$}

Since there is a lattice structure on packed words and since we know that
the product $\M_u \M_v$ is an interval of this lattice, we can 
define several interesting bases, depending on the way we use the lattice.

As in the case of the permutohedron, one can take sums of $\M_u$, over all the
elements upper or lower than $u$ in the lattice, or restricted  to 
elements belonging to the same ``class'' as $u$ (see~\cite{NCSF6,AS} for 
examples of such bases). In the case of the
permutohedron, the classes are the descent classes of permutations. 
In our case, the classes
are the intervals of the pseudo-permutohedron composed of words with the same
standardization.

Summing over all elements upper (or lower) than a word $u$ naturally yields
multiplicative bases on $\WQSym$. Summing over all elements upper (or
lower) than $u$ inside its standardization class leads to  analogs of the
usual bases of $\QSym$.

\medskip
%%%%%%%%%%%%%%%%%%%%%%%%%%%%%%%%%%%%
\subsubsection{Multiplicative bases}

Let
\begin{equation}
\SW_u := \sum_{v\pluspetit u} \M_v \quad\text{and}\quad
\EW_u := \sum_{u\pluspetit v} \M_v.
\end{equation}

For example,
\begin{equation}
\SW_{212} = \M_{212} + \M_{213} + \M_{112} + \M_{123}.
\end{equation}
\begin{equation}
\EW_{212} = \M_{212} + \M_{312} + \M_{211} + \M_{321}.
\end{equation}
\begin{equation}
\SW_{1122} = \M_{1122} + \M_{1123} + \M_{1233} + \M_{1234}.
\end{equation}

Since both $\SW$ and $\EW$ are triangular over the basis $\M_u$ of $\WQSym$,
we know that these are bases of $\WQSym$.

\begin{theorem}
The sets $(\SW_u)$ and $(\EW_u)$ where $u$ runs over packed words
are bases of $\WQSym$. Moreover, their product is given by

\begin{equation}
\SW_{u'} \SW_{u''} = \SW_{u'[\max(u'')]\cdot u''}.
\end{equation}
\begin{equation}
\EW_{u'} \EW_{u''} = \EW_{u'\cdot u''[\max(u')]}.
\end{equation}
\end{theorem}

For example,
\begin{equation}
%(wqs::G([1,1,2,2]) + wqs::G([1,1,2,3]) + wqs::G([1,2,3,3]) +
%wqs::G([1,2,3,4])) * (wqs::G([1,3,2]) + wqs::G([1,2,2]) + wqs::G([1,2,3]));
%-> 4455132
\SW_{1122}\SW_{132} = \SW_{4455132}.
\end{equation}
\begin{equation}
\EW_{1122}\EW_{132} = \EW_{1122354}.
\end{equation}

\medskip
%%%%%%%%%%%%%%%%%%%%%%%%%%%%%%%%%%%%%%%%%%%%%%
\subsubsection{Quasi-ribbon basis of $\WQSym$}

Let us first mention that a basis of $\WQSym$ has been defined in~\cite{BZ}
by summing over intervals restricted to standardization classes of packed words. 

We will now consider similar sums but taken the other way round, in order to
build the analogs of $\WQSym$ of Gessel's fundamental basis $F_I$ of
$\QSym$. Indeed, as already
mentioned, the $\M_u$ are mapped to the $M_I$ of $\QSym$ under the
abelianization $\K\langle A\rangle\rightarrow\K[X]$ of $\WQSym$. 
Since the pair of dual bases $(F_I,R_I)$ of $(\QSym,\NCSF)$ 
is of fundamental importance, it is natural to ask whether one can find an
analogous pair for $(\WQSym,\WQSym^*)$. To avoid confusion in the
notations, we will denote the analog of  $F_I$  by $\Phi_u$ instead of
$\F_u$ since this notation is already used in the dual algebra
$\WQSym^*\subset \PQSym$, with a different meaning.
The analog of $R$ basis in $\WQSym^*$ will still be denoted by $R$.
The representation of packed words by segmented permutations is more suited
for the next statements since one easily checks that two words $u$ and $v$
having the same standardized word satisfy $v\pluspetit u$ iff $v$ is obtained
as a segmented permutation from the segmented permutation of $u$ by inserting
any number of bars.
Let
%Recall that the \emph{reverse refinement order}, denoted by $\raff$, on
%compositions is such that $I=(i_1,\ldots,i_k)\raff J=(j_1,\ldots,j_l)$ iff
%$\{i_1,i_1+i_2,\ldots,i_1+\cdots+i_k\}$ contains
%$\{j_1,j_1+j_2,\ldots,j_1+\cdots+j_l\}$.
%In this case, we say that $I$ is finer than $J$.
%For example, $(2,1,2,3,1,2)\raff (3,2,6)$.
%
%This  order can be extended to packed words as follows.
%We say that $w$ is finer than $w'$, and  write $w\raff w'$, iff
%$w$ and $w'$ have same standardized word and the evaluation of $w$ is finer
%than the evaluation of $w'$.
%
%Let
%\begin{equation} 
%\Phi_u := \sum_{v ; v\raff u} \M_v.
%\end{equation} 
%On the representation of packed words by segmented permutations, let
\begin{equation} 
\Phi_\sigma := \sum_{\sigma'} \M_{\sigma'}
\end{equation}
where $\sigma'$ runs ver the set of segmented permutations obtained from
$\sigma$ by inserting any number of bars.
For example,
%\begin{equation}
%\Phi_{111} = \M_{111} + \M_{112} + \M_{122} + \M_{123} ; \quad
%\Phi_{112} = \M_{112} + \M_{123} ;
%\Phi_{212} = \M_{212} + \M_{213}.
%\end{equation}
%\begin{equation}
%\Phi_{133142} = \M_{133142} + \M_{134152} + \M_{144253} + \M_{145263}.
%\end{equation}
\begin{equation}
\Phi_{14|6|23|5} = \M_{14|6|23|5} + \M_{14|6|2|3|5} + \M_{1|4|6|23|5} +
\M_{1|4|6|2|3|5}.
\end{equation}

Since $(\Phi_u)$ is triangular over  $(\M_u)$, it is a basis of
$\WQSym$.
% Note that the order used for summation is a restriction
%of the refinement order on compositions, so is a boolean lattice. 
%Hence,
%%
%\begin{equation}
%\M_u = \sum_{v; v\raff u} (-1)^{\max(v)-\max(u)} \Phi_v.
%\end{equation}
%%
%For example,
%\begin{equation}
%\M_{133142} = \Phi_{133142} - \Phi_{134152} - \Phi_{144253} + \Phi_{145263}.
%\end{equation}
%
By construction, it satisfies a product formula  similar to
that of Gessel's basis $F_I$ of $\QSym$ (whence the choice of notation).
To state it, we need an analogue of the shifted shuffle, defined on the special
class of segmented permutations encoding set compositions.

The \emph{shifted shuffle} $\alpha\ssh\beta$ of two such segmented permutations is
obtained from the usual shifted shuffle $\sigma\ssh\tau$ of the underlying
permutations $\sigma$ and $\tau$ by inserting bars
\begin{itemize}
\item between each pairs of letters coming from the same word if they were
separated by a bar in this word,
\item after each element of $\beta$ followed by an element of $\alpha$.
\end{itemize}

For example,
\begin{equation}
2|1 \ssh 12 = 2|134 + 23|14 + 234|1 + 3|2|14 + 3|24|1 + 34|2|1.
\end{equation}
%
%We then have:
%
\begin{theorem}
The product and coproduct in the basis $\Phi$ are given by
\begin{equation}
\Phi_{\sigma'} \Phi_{\sigma''} = \sum_{\sigma\in \sigma'\ssh\sigma''}
\Phi_{\sigma}.
\end{equation}
\begin{equation}
\Delta\Phi_\sigma =
\sum_{\sigma'|\sigma''=\sigma \text{\ or\ } \sigma'\cdot  \sigma''=\sigma}
\Phi_{\Std(\sigma')} \otimes \Phi_{\Std(\sigma'')}.
\end{equation}
\end{theorem}

%For example, in both encodings, we have
For example, we have
\begin{equation}
%\begin{split}
%& \Phi_{1} \Phi_{121} = \Phi_{1121} + \Phi_{2132} + \Phi_{2121} + \Phi_{3121}.
%\\
%&
\Phi_{1} \Phi_{13|2} = \Phi_{124|3} + \Phi_{2|14|3} + \Phi_{24|13} +
\Phi_{24|3|1}.
%\end{split}
\end{equation}
%\begin{equation}
%\begin{split}
%wqs::Phi(wqs::Phi([1,3,1,2]) * wqs::Phi([2,1]));
%\Phi_{1312} \Phi_{21} =& \ \ \ \
% \Phi_{131221} + \Phi_{131231} + \Phi_{131232} + \Phi_{131243} +
% \Phi_{141232} \\ &+
% \Phi_{141321} + \Phi_{142321} + \Phi_{142331} + \Phi_{142341} +
% \Phi_{153421} \\ &+
% \Phi_{242321} + \Phi_{242331} + \Phi_{242341} + \Phi_{253421} +
% \Phi_{353421}.
%\end{split}
%\end{equation}
%
%\begin{equation}
%\Delta\Phi_{23121} =
%1\otimes\Phi_{23121} + \Phi_{1}\otimes\Phi_{2321} +
%\Phi_{11}\otimes\Phi_{121} + \Phi_{211}\otimes\Phi_{21} +
%\Phi_{2121}\otimes\Phi_{1} + \Phi_{23121}\otimes1.
%\end{equation}
%
\begin{equation}
\Delta\Phi_{35|14|2} =
1\otimes\Phi_{35|14|2} + \Phi_{1}\otimes\Phi_{4|13|2} +
\Phi_{12}\otimes\Phi_{13|2} + \Phi_{23|1}\otimes\Phi_{2|1} +
\Phi_{24|13}\otimes\Phi_{1} + \Phi_{35|14|2}\otimes1.
\end{equation}

Note that under abelianization, $\chi(\Phi_u)=F_I$ where $I$ is the evaluation
of $u$.
%where $\Phi$ is sent to the basis $F$ of $\QSym$.

\medskip
%%%%%%%%%%%%%%%%%%%%%%%%%%%%%%%%%%%%%%%%%%
\subsubsection{Ribbon basis of $\WQSym^*$}

Let us now consider the dual basis of $\Phi$. We have seen that
it should be regarded as an analog of the ribbon basis of $\Sym$.
%Even if the product of our analogs of the $R_I$ contains more terms than the
%usual $R_I$, we will see later that they build up the right ribbon functions
%in subalgebras and quotients of $\WQSym^*$.
%
By duality, one can state:

\begin{theorem}
Let $R_\sigma$ be the dual basis of $\Phi_\sigma$. Then the product and
coproduct in this basis are given by
\begin{equation}
R_{\sigma'} R_{\sigma''} = \sum_{\sigma=\tau|\nu \text{\ or\ } \sigma=\tau\nu;
\Std(\tau)=\sigma', \Std(\nu)=\sigma''} R_{\sigma}.
\end{equation}
\begin{equation}
\Delta R_\sigma = \sum_{\sigma'.\sigma''=\sigma}
R_{\Std(\sigma')} \otimes R_{\Std(\sigma'')}.
\end{equation}
\end{theorem}

Note that there are more elements coming from $\tau|\nu$ than from $\tau\nu$
since the permutation $\sigma$ has to be \emph{increasing} between two bars.

For example,
\begin{equation}
%wqs::Ra(wqs::Ra([2,1]) * wqs::Ra([1]));    
R_{21} R_{1} = R_{212} + R_{221} + R_{213} + R_{231} + R_{321}.
\end{equation}
%\begin{equation}
%%wqs::Ra(wqs::Ra([1,2,1]) * wqs::Ra([1,2]));
%\begin{split}
%R_{121} R_{12} =& \ \ \ \
%R_{12123} + R_{12213} + R_{12231} + R_{12134} + R_{12314} \\ &+
%R_{12341} + R_{13214} + R_{13241} + R_{13421} + R_{31214} \\ &+
%R_{31241} + R_{31421} + R_{34121}.
%\end{split}
%\end{equation}
%\begin{equation}
%\begin{split}
%%wqs::Ra(wqs::Ra([1,2,1]) * wqs::Ra([2,1]));
%R_{121} R_{21} =&\ \ \ \
%R_{12132} + R_{12312} + R_{12321} + R_{13212} + R_{13221} + R_{31212} \\ &+
%R_{31221} + R_{12143} + R_{12413} + R_{12431} + R_{14213} + R_{14231} \\ &+
%R_{14321} + R_{41213} + R_{41231} + R_{41321} + R_{43121}.
%\end{split}
%\end{equation}

%%%%%%%%%%%%%%%%%%%%%%%%%%%%%%%%%%%%%%%%%%%%%%%%%%%%%%%%%%%%%%%%%%%%%%%%%%%%%%%
%%%%%%%%%%%%%%%%%%%%%%%%%%%%%%%%%%%%%%%%%%%%%%%%%%%%%%%%%%%%%%%%%%%%%%%%%%%%%%%
%%%%%%%%%%%%%%%%%%%%%%%%%%%%%%%%%%%%%%%%%%%%%%%%%%%%%%%%%%%%%%%%%%%%%%%%%%%%%%%
\section{Hopf algebras based on Schr\"oder sets} 
\label{td}
\label{sqsym}

In Section~\ref{trigd1}, we recalled that the little Schr\"oder numbers build
up the Hilbert series of the free dendriform trialgebra on one generator $\TD$.
We will see that our relization of $\TD$ endows it with a natural structure of
bidendriform bialgebra.  In particular, this will prove that there is a
natural self-dual Hopf structure on $\TD$.
But there are other ways to arrive at the little Schr\"oder numbers from the
other Hopf algebras $\WQSym$ and $\PQSym$. Indeed, the number of classes of
packed words of size $n$ under the sylvester congruence is $s_n$, and the
number of classes of parking functions of size $n$ under the hypoplactic
congruence is also $s_n$.
The hypoplactic quotient of $\PQSym^*$ has been studied in~\cite{NT3}. It is
not isomorphic to $\TD$ nor to the sylvester quotient of $\WQSym$ since it is
a non self-dual Hopf algebra whereas the last two are self-dual, and
furthemore isomorphic as bidendriform bialgebras and as dendriform
trialgebras.

%%%%%%%%%%%%%%%%%%%%%%%%%%%%%%%%%%%%%%%%%%%%%%%%%%%%%%%%%%%%%%%%%%%%%%%%%%%%%%%
\subsection{The free dendriform trialgebra again}

Recall that we realized the free dendriform trialgebra in
Section~\ref{trigd1} as the subtrialgebra of
$\WQSym$ generated by $\M_1$, the sum of all letters.
It is immediate that $\TD$ is stable by the codendriform
half-coproducts of $\WQSym^*$. Hence,

\begin{theorem}
\label{free-all}
$\TD$  is a sub-bidendriform bialgebra, and hence
a Hopf subalgebra of $\WQSym^*$.
In particular, $\TD$ is free, self-dual and its primitive Lie algebra is
free.
\end{theorem}

%%%%%%%%%%%%%%%%%%%%%%%%%%%%%%%%%%%%%%%%%%%%%%%%%%%%%%%%%%%%%%%%%%%%%%%%%%%%%%%
\subsection{Lattice structure on plane trees}

Given a plane tree $T$, define its \emph{canonical word} as the
maximal packed word $w$ in the pseudo-permutohedron such that $\TT(w)=T$.

%For example, the canonical words up to $n=4$ are
For example, the canonical words up to $n=3$ are
\begin{equation}
\label{ex1234}
\begin{split}
&\{1\},\qquad \{11, 12, 21\},\qquad
\{111,\ 112,\ 211,\ 122,\ 212,\ 221,\ 123,\ 213,\ 231,\ 312,\ 321\} \\
%&\{1111,\ 1112,\ 2111,\ 1122,\ 2112,\ 2211,\ 1123,\ 2113,\ 2311,\ 3112,\
%3211, 1222,\\ &
%2122,\ 2212,\ 2221,\ 1223,\ 2123,\ 2213,\ 2231,\ 3122,\ 3212,\ 3221, 1233,\
%2133,\\ &
%2313,\ 2331,\ 3123,\ 3213,\ 3231,\ 3312,\ 3321,\ 1234,\ 2134,\ 2314,\ 2341,\
%3124,\\ & 3214,\ 3241,\ 3412,\ 3421,\ 4123,\ 4213,\ 4231,\ 4312,\ 4321\}.
\end{split}
\end{equation}
%The lattice for $n=3$ is represented on Figure~\ref{lat-pt}.

\begin{figure}[ht]
\resizebox{3cm}{4cm}{\includegraphics{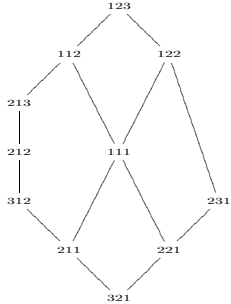}}
\caption{\label{lat-pt}The lattice of plane trees represented by their
canonical words for $n=3$.}
\end{figure}

%Starting from a tree $T$, one gets its canonical word as follows. Let
%$T_1,\ldots,T_k$ be the subtrees of $T$  of the root
%of $T$. Let $w_i$ be the canonical word of each $T_i$. Then
%\begin{equation}
%w = w_1[\max(w_2)+\cdots+\max(w_k)]\, m\, w_2[\max(w_3)+\cdots+\max(w_k)]\,
%m \ldots m\, w_{k-1}[\max(w_k)]\, m\, w_k,
%\end{equation}
%where $m=\max(w_1)+\cdots+\max(w_k)$.

Define the \emph{second canonical word} of each tree $T$ as the
minimal packed word $w$ in the pseudo-permutohedron such that $\TT(w)=T$.
%There is a similar algorithm to compute this word from a tree.
% second canonical word is computed in a similar way.

A packed word $u=u_1\cdots u_n$ is said to \emph{avoid the pattern}
$w=w_1\cdots w_k$ if there is no sequence
$1\leq i_1<\cdots<i_k\leq n$ such that $u' = u_{i_1}\cdots u_{i_k}$ has same
inversions and same half-inversions as $w$.

For example, $41352312$ avoids the patterns $111$ and $1122$, but not $2311$
since $3522$ has the same (half)-inversions.

\begin{theorem}
The canonical words of trees are the packed words avoiding the patterns 121
and 132.
The second canonical words of trees are the packed words avoiding the patterns
121 and 231.
\end{theorem}

Set $u\sim_T v$ iff $\TT(u)=\TT(v)$.
We now define two orders $\sim_T$-classes of packed words
\begin{itemize}
\item[1.] A class $S$ is smaller than a class $S'$ if the canonical word of
$S$ is smaller than the canonical word of $S'$ in the pseudo-permutohedron.
\item[2.] A class $S$ is smaller than a class $S'$ if there is a pair $(w,w')$
in $S\times S'$ such that $w$ is smaller than $w'$ in the
pseudo-permutohedron.
\end{itemize}

\begin{theorem}
These two orders coincide and are also identical with the one defined
in~\cite{PR}. Moreover, the restriction of the pseudo-permutohedron to the
canonical words of trees is a lattice.
\end{theorem}

%%%%%%%%%%%%%%%%%%%%%%%%%%%%%%%%%%%%%%%%%%%%%%%%%%%%%%%%%%%%%%%%%%%%%%%%%%%%%%%
\subsection{Some bases of $\TD$}

\subsubsection{The basis $\MM_T$}

Let us start with the already defined basis $\MM_T$.
First note that $\MM_{T}$ expressed as a sum of $\M_u$ in $\WQSym$ is
an interval of the pseudo-permutohedron. From the above description of the
lattice, we obtain easily:

\begin{theorem}[\cite{PR}]
The product $\MM_{T'}\MM_{T''}$ is an interval of the  lattice of plane trees.
On trees, the minimum $T'\wedge T''$ is obtained by
gluing the root of $T''$ at the end of the
leftmost branch of $T'$, whereas the maximum $T'\vee T''$
is obtained by gluing the root of
$T'$ at the end of the rightmost branch of $T''$.

On the canonical words $w'$ and $w''$, the minimum is the canonical word
associated with $w'\cdot w''[\max(w')]$ and the maximum is
$w'[\max(w'')]\cdot w''$.
\end{theorem}

%%%%%%%%%%%%%%%%%%%%%%%%%%%%%%%%%%%%%%%%%%%%%%%%%%%%%%%%%%%%%%%%%%%%%%%%%%%%%%%
\subsubsection{Complete and elementary bases of $\TD$}

We can also build two multiplicative bases as  in $\WQSym$.
% They are particularly simple:

\begin{theorem}
The set $(\SW_{w})$ (resp. $(\EW_{w})$) where $w$ runs over canonical (resp.
second canonical) words are multiplicative bases of $\TD$.
%The set $(\EW_{w})$ where $w$ runs over second canonical words is a
%basis of $\TD$.
\end{theorem}

%%%%%%%%%%%%%%%%%%%%%%%%%%%%%%%%%%%%%%%%%%%%%%%%%%%%%%%%%%%%%%%%%%%%%%%%%%%%%%%
\subsection{Internal product on $\TD$}

If one defines $\TD$ as the Hopf subalgebra of $\WQSym$ defined by
\begin{equation}
\MM_T=\sum_{\TT(u)=T}\M_u\,,
\end{equation}
then $\TD^*$ is the quotient of $\WQSym^*$ by the relation $\F_u\equiv\F_v$
iff $\TT(u)=\TT(v)$. We denote by $S_T$ the dual basis of $\MM_T$.

\begin{theorem}
The internal product of $\WQSym^*_n$ induces an internal product on the
homogeneous components $\TD^*_n$ of the dual  algebra. 
More precisely, one has
\begin{equation}
\label{prodintS}
S_{T'} * S_{T''} = S_T,
\end{equation}
where $T$ is the tree obtained by applying $\TT$ to the biword of the
canonical words of the trees $T'$ and $T''$.
\end{theorem}

For example, representing trees as their canonical words, one has
\begin{equation}
S_{221} * S_{122} = S_{231}; \qquad S_{221}*S_{321} = S_{321};
\end{equation}
%\begin{equation}
%S_{113} * S_{231} = S_{123}; \qquad S_{311}*S_{311} = S_{311};
%\end{equation}
\begin{equation}
S_{453223515} * S_{433442214} = S_{674223518}.
\end{equation}

%%%%%%%%%%%%%%%%%%%%%%%%%%%%%%%%%%%%%%%%%%%%%%%%%%%%%%%%%%%%%%%%%%%%%%%%%%%%%%%
\subsection{Sylvester quotient of $\WQSym$}

One can check by direct calculation that the sylvester quotient~\cite{HNT} of
$\WQSym$
is also stable by the tridendriform operations, and by the codendriform
half-coproducts since the elements of a sylvester class have the same last
letter. Hence,

\begin{theorem}
The sylvester quotient of $\WQSym$ is a dendriform trialgebra, a
bidendriform bialgebra, and hence a Hopf algebra.
It is isomorphic to $\TD$ as a dendriform trialgebra, as a bidendriform
bialgebra and as a Hopf algebra.
\end{theorem}

%%%%%%%%%%%%%%%%%%%%%%%%%%%%%%%%%%%%%%%%%%%%%%%%%%%%%%%%%%%%%%%%%%%%%%%%%%%%%%%
%%%%%%%%%%%%%%%%%%%%%%%%%%%%%%%%%%%%%%%%%%%%%%%%%%%%%%%%%%%%%%%%%%%%%%%%%%%%%%%
%%%%%%%%%%%%%%%%%%%%%%%%%%%%%%%%%%%%%%%%%%%%%%%%%%%%%%%%%%%%%%%%%%%%%%%%%%%%%%%
\section{A Hopf algebra of segmented compositions}
\label{cc}

In~\cite{NT3}, we have built a Hopf subalgebra $\SCQSym^*$ of the hypoplactic
quotient $\SQSym^*$ of $\PQSym^*$, whose Hilbert series is given by
\begin{equation}
1 + \sum_{n\geq1} 3^{n-1}t^n.
\end{equation}
This Hopf algebra is not self-dual, but admits lifts of Gessel's
fundamental basis
$F_I$ of $\QSym$ and its dual basis.
Since the elements of $\SCQSym^*$ are obtained by summing up hypoplactic
classes having the same packed word, thanks to the following diagram,
it is obvious that $\SCQSym^*$ is also the quotient of $\WQSym$ by the
hypoplactic congruence.

\begin{equation}
\begin{CD}
\PQSym^* @>hypo>> \SQSym^* \\
@A(pack)AA @AA(pack)A \\
\WQSym   @>hypo>> \SCQSym^*
\end{CD}
\end{equation}

%%%%%%%%%%%%%% Délire à partir d'ici %%%%%%%%%%%%%%%%%%%%%%%%%%
%{\tt D\'elirant ?}
%However, it is possible to obtain another Hopf algebra $\CC$ indexed by the
%same objects that is a bidendriform bialgebra, thus self-dual as a Hopf
%algebra.
%As in the case of $\TD$, this algebra is related to a lattice that will be a
%sublattice a quotient lattice of the lattice built in Section~\ref{td} and
%of the pseudo-permutohedron. In particular, $\CC$ is a Hopf subalgebra of
%$\TD$ with multiplicative bases built as before.

%%%%%%%%%%%%%%%%%%%%%%%%%%%%%%%%%%%%%%%%%%%%%%%%%%%%%%%%%%%%%%%%%%%%%%%%%%%%%%%
\subsection{Segmented compositions}

Define a \emph{segmented composition} as a finite sequence of integers,
separated by vertical bars or commas, \emph{e.g.}, $(2,1\sep2\sep1,2)$.

The number of segmented compositions having the same underlying composition is
obviously $2^{l-1}$ where $l$ is the length of the composition, so that the
total number of segmented compositions of sum $n$ is $3^{n-1}$.
There is a natural bijection between segmented
compositions of $n$ and sequences of length $n-1$ over three symbols
$<,=,>$: start with a segmented composition $\Ig$. If the $i$-th position is
not a descent of the underlying ribbon diagram, write $<$ ; otherwise, if $i$
is followed by a comma, write $=$ ; if $i$ is followed by a bar, write $>$.

Now, with each word $w$ of length $n$,  associate a segmented
composition $\SC(w)=s_1 \cdots s_{n-1}$ where $s_i$ is the
correct comparison sign between $w_i$ and $w_{i+1}$.
For example, given $w=1615116244543$, one gets the sequence (and the segmented
composition):
\begin{equation}
<><>=<><=<>>  \Longleftrightarrow (2|2|1,2|2,2|1|1).
\end{equation}

%%%%%%%%%%%%%%%%%%%%%%%%%%%%%%%%%%%%%%%%%%%%%%%%%%%%%%%%%%%%%%%%%%%%%%%%%%%%%%%
\subsection{A subalgebra of $\TD$}

Given a segmented composition $\Ig$, define
\begin{equation}
\label{wt2cc}
\MC_\Ig=\sum_{\SC(T)=\Ig}\MM_T = \sum_{\SC(u)=\Ig}\M_u\,.
\end{equation}

For example,
\begin{equation}
%select(listCanonsDend(4), z-> (z[1]=z[2] and z[2]<z[3] and z[3]>z[4]));
%select(listCanonsDend(4), z-> (z[1]>z[2] and z[2]<z[3] and z[3]<z[4]));
\MC_{12|1} = \MM_{2231} \qquad
\MC_{1|3} = \MM_{2123} + \MM_{2134} + \MM_{3123} + \MM_{3124} + \MM_{4123}.
\end{equation}
%\begin{equation}
%%select(listCanonsDend(5), z-> (z[1]>z[2] and z[2]<z[3] and z[3]>z[4] and
%%z[4]=z[5]));
%\MC_{1|2|11} = \MM_{32311} + \MM_{32411} + \MM_{42311}.
%\end{equation}

%On the realization on words, one sees that
%\begin{equation}
%\MC_\Ig=\sum_{\SC(w)=\Ig}w,
%\end{equation}
%so that 
%
\begin{theorem}
The $\MC_\Ig$ generate a subalgebra $\TC$ of $\TD$.
Their product is given by
\begin{equation}
\MC_{\Ig'} \MC_{\Ig''} = \MC_{\Ig'\triangleright \Ig''}+\MC_{\Ig',\Ig''}+\MC_{\Ig'|\Ig''}.
\end{equation}
where $\Ig'\triangleright \Ig''$ is obtained by gluing  the last part of
$\Ig'$ and the first part of $\Ig''$, so that $\TC$ is the free cubical
trialgebra on one generator~\cite{LRtri}.
\end{theorem}

For example,
\begin{equation}
\MC_{1|21}\MC_{31} = \MC_{1|241} + \MC_{1|2131} + \MC_{1|21|31}.
\end{equation}

%Moreover, $\CC$ is stable by the codendriform half-coproducts of $\TD$.
%Hence,

%\begin{theorem} 
%\label{freec-all}
%$\CC$ is a sub-dendriform trialgebra, a sub-bidendriform bialgebra, and hence
%a Hopf subalgebra of $\TD$.
%In particular, $\CC$ is free, self-dual and its primitive Lie algebra is
%primitive.
%\end{theorem}

%%%%%%%%%%%%%%%%%%%%%%%%%%%%%%%%%%%%%%%%%%%%%%%%%%%%%%%%%%%%%%%%%%%%%%%%%%%%%%%
\subsection{A lattice structure on segmented compositions}

Given a segmented composition $\Ig$, define its \emph{canonical word} as the
maximal packed word $w$ in the pseudo-permutohedron such that $\SC(w)=\Ig$.

For example, the canonical words up to $n=3$ are
\begin{equation}
\label{exc1234}
\begin{split}
&\{1\},\qquad \{11, 12, 21\},\qquad
\{111,\ 112,\ 211,\ 122,\ 221,\ 123,\ 231,\ 312,\ 321\} \\
%&\{1111,\ 1112,\ 2111,\ 1122,\ 2211,\ 1123,\ 2311,\ 3112,\ 3211,\\ &
%   1222,\ 2221,\ 1223,\ 2231,\ 3122,\ 3221,\ 1233,\ 2331,\ 3312,\\ &
%   3321,\ 1234,\ 2341,\ 3412,\ 3421,\ 4123,\ 4231,\ 4312,\ 4321\}.
\end{split}
\end{equation} 
%The lattice for $n=3$ is represented on Figure~\ref{lat-sc}.

\begin{figure}[ht]
\resizebox{3cm}{4cm}{\includegraphics{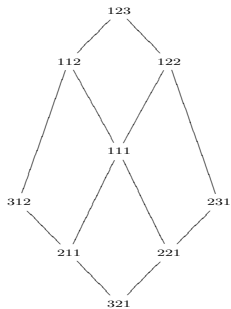}}
\caption{\label{lat-sc}The lattice of segmented compositions represented by
their canonical words at $n=3$.}
\end{figure}

%Starting from a segmented composition $\Ig$ of $n$
%written as a sequence $s$ of $<,=,>$,
%one gets its canonical word as follows. Let $i=n-|s|_=$,
%where $|s|_=$ is the number of symbols $=$ in $s$.
%Put $i$ at the
%leftmost position such that it can be the maximum of the canonical word
%(inside a sequence $<==\cdots=>$). Then put $i-1$ at the leftmost position
%that does not contradict the rule, and so on.
%For example, with the sequence $=<=>>=<=$, one finds $i=5$ and place the
%$5$'s as $=<5=5>>=<=$, then the $4$'s: $4=4<5=5>>=<=$, then $3$, $2$, and
%$1$, to get $445531122$.

Define the \emph{second canonical word} of a segmented composition $\Ig$ as
the minimal packed word $w$ in the pseudo-permutohedron such that $\SC(w)=\Ig$.
%There is a similar algorithm to compute this word from a segmented
%composition.

\begin{theorem}
The canonical words of segmented compositions are the packed words avoiding
the patterns 121, 132, 212, and 213.
The second canonical words of segmented compositions are the packed words
avoiding the patterns 121, 231, 212, and 312.
\end{theorem}

Let $u\sim_S v$ iff $\SC(u)=\SC(v)$.
We  define two orders on $\sim_S$-equivalence classes of words. 

\begin{itemize}
\item[1.] A class $S$ is smaller than a class $S'$
if the canonical word of $S$ is smaller than the canonical word of $S'$ in the
pseudo-permutohedron.
\item[2.] A class $S$ is smaller than a class $S'$
if there exists two elements $(w,w')$ in $S\times S'$ such that
$w$ is smaller than $w'$ in the pseudo-permutohedron.
\end{itemize}

\begin{proposition}
The two orders coincide. Moreover, the restriction of the pseudo-permutohedron
to the canonical segmented words is a lattice.
\end{proposition}

%%%%%%%%%%%%%%%%%%%%%%%%%%%%%%%%%%%%%%%%%%%%%%%%%%%%%%%%%%%%%%%%%%%%%%%%%%%%%%%
\subsection{Multiplicative bases}

We can  build two multiplicative bases, as  in $\WQSym$. They are
particularly simple:

\begin{theorem}
The set $(\SW_{w})$ where $w$ runs into the set of canonical segmented words
is a basis of $\TC$.
The set $(\EW_{w})$ where $w$ runs into the set of second canonical segmented
words is a basis of $\TC$.
\end{theorem}

%%%%%%%%%%%%%%%%%%%%%%%%%%%%%%%%%%%%%%%%%%%%%%%%%%%%%%%%%%%%%%%%%%%%%%%%%%%%%%%
\subsection{Internal product on $\TC$}

%If one defines $\TC$ as the Hopf subalgebra of $\WQSym$ defined by
If one defines $\TC$ as the Hopf subalgebra of $\WQSym$ as in
Equation~(\ref{wt2cc}),
%\begin{equation}
%\MM_\Ig=\sum_{\SC(u)=\Ig}\M_u\,,
%\end{equation}
then $\TC^*$ is the quotient of $\WQSym^*$ by the relation $\F_u\equiv \F_v$ iff
$\SC(u)=\SC(v)$. We denote by $S_\Ig$ the dual basis of $\MC_\Ig$.

\begin{theorem}
The internal product of $\WQSym^*$ induces an internal product on
the homogeneous components $\TC^*_n$ of $\TC^*$. 
More precisely, one has
\begin{equation}
\label{prodintcS}
S_{\Ig'} * S_{\Ig''} = S_\Ig,
\end{equation}
where $\Ig$ is the segmented composition obtained by applying $\SC$ to the
biword of the canonical words of the segmented compositions $\Ig'$ and
$\Ig''$.
\end{theorem}

%%%%%%%%%%%%%%%%%%%%%%%%%%%%%%%%%%%%%%%%%%%%%%%%%%%%%%%%%%%%%%%%%%%%%%%%%%%%%%%
%%%%%%%%%%%%%%%%%%%%%%%%%%%%%%%%%%%%%%%%%%%%%%%%%%%%%%%%%%%%%%%%%%%%%%%%%%%%%%%
%%%%%%%%%%%%%%%%%%%%%%%%%%%%%%%%%%%%%%%%%%%%%%%%%%%%%%%%%%%%%%%%%%%%%%%%%%%%%%%
%%%%%%%%%%%%%%%%%%%%%  BIBLIOGRAPHIE !!! %%%%%%%%%%%%%%%%%%%%%%%%%%%%%%%%%%%%%%
%%%%%%%%%%%%%%%%%%%%%%%%%%%%%%%%%%%%%%%%%%%%%%%%%%%%%%%%%%%%%%%%%%%%%%%%%%%%%%%
%%%%%%%%%%%%%%%%%%%%%%%%%%%%%%%%%%%%%%%%%%%%%%%%%%%%%%%%%%%%%%%%%%%%%%%%%%%%%%%

\end{document}